\documentclass[10pt]{article}
\usepackage[cp866]{inputenc}
\begin{document}
\bigskip 

\centerline{{\Large\bf Role of skew-symmetric differential 
forms }}
\centerline{{\Large\bf in mathematics }}
\centerline {\bf L.~I. Petrova}
\centerline{{\it Moscow State University, Russia, e-mail: ptr@cmc.msu.ru}}

\renewcommand{\abstractname}{Abstract}
\begin{abstract}

Skew-symmetric forms possess unique capabilities. This is due to the fact
that they deal with differentials
and differential expressions; and, therefore, they are suitable for
describing invariants and invariant structures. The properties of closed 
exterior and dual forms, namely, invariance, covariance, conjugacy and duality, 
either explicitly or implicitly appear in all invariant mathematical formalisms.
The closed exterior forms relate to such branches of mathematics as
algebra, geometry, differential geometry, vectorial and tensor calculus,
the theory of complex variables, differential equations and so on.
This enables one to see an internal connection between various branches
of mathematics.

However, the theory of closed exterior forms cannot be completed
without an answer to a question of how the closed exterior forms
emerge. In the present paper we discus essentially new
skew-symmetric forms, which possess capabilities so far not
exploited in any mathematical formalisms, namely, they generate
closed exterior forms.
Such skew-symmetric forms, which are evolutionary ones, are derived from
differential equations, and, in contrast to exterior forms, they are defined
on nonintegrable manifolds. Theory of evolutionary forms, what includes such 
elements as nonidentical relations, degenerate transformations, 
transition from nonintegrable manifold to integrable one,  
enables one to understand the process of conjugating operators,
the mechanism of generation of invariant structures and others. 

\end{abstract}

\section{Exterior differential forms}

In this Section introductory information on exterior differential  forms is
presented, and the basic properties of closed exterior  differential forms
and specific features of their mathematical apparatus are described. 

The concept of ``Exterior differential forms  (differential forms with
exterior multiplication)" was introduced by E.Cartan as a general concept for
integrands, which constitute integral invariants [1]. (The existence of integral
invariants was recognized previously by A. Poincare while studying the
general equations of dynamics.)

The exterior differential form of degree $p$ ($p$-form) can be
written as [2--4]
$$
\theta^p=\sum_{i_1\dots i_p}a_{i_1\dots i_p}dx^{i_1}\wedge
dx^{i_2}\wedge\dots \wedge dx^{i_p}\quad 0\leq p\leq n\eqno(1.1)
$$
Here $a_{i_1\dots i_p}$ are the functions of the variables $x^{i_1}$,
$x^{i_2}$, \dots, $x^{i_n}$, $n$ is the dimension of the space,
$\wedge$ is the operator of exterior multiplication, $dx^i$,
$dx^{i}\wedge dx^{j}$, $dx^{i}\wedge dx^{j}\wedge dx^{k}$, \dots\
is the local basis which satisfies the condition of exterior
multiplication:
$$
\begin{array}{l}
dx^{i}\wedge dx^{i}=0\\
dx^{i}\wedge dx^{j}=-dx^{j}\wedge dx^{i}\quad i\ne j
\end{array}\eqno(1.2)
$$
[Below, the symbol for summing, $\sum$, and the symbol for exterior
multiplication, $\wedge$, will be omitted. Summation over repeated indices is
implied.]

The differential of the exterior form $\theta^p$ is expressed as
$$
d\theta^p=\sum_{i_1\dots i_p}da_{i_1\dots
i_p}dx^{i_1}dx^{i_2}\dots dx^{i_p} \eqno(1.3)
$$
and is a differential form of degree $(p+1)$. 

Let us consider some examples of the exterior differential form whose
basis is the Euclidean space domains.

We consider a 3-dimensional space. In this case the differential
forms of zero-, first- and second degree can be written as, [4]:
$$\theta^0=a$$
$$\theta^1=a_1 dx^1+a_2 dx^2+a_3 dx^3$$
$$\theta^2=a_{12}dx^1 dx^2+a_{23}dx^2 dx^3+a_{31}dx^3 dx^1$$

With account for conditions (1.2), their differentials are the forms

$$d\theta^0=\frac{\partial a}{\partial x^1}dx^1+\frac{\partial a}{\partial
x^2}dx^2+\frac{\partial a}{\partial x^3}dx^3$$
$$d\theta^1=\left(\frac{\partial a_2}{\partial x^1}-\frac{\partial
a_1}{\partial x^2}\right)dx^1 dx^2+\left(\frac{\partial a_3}{\partial
x^2}-\frac{\partial a_2}{\partial x^3}\right)dx^2 dx^3+\nonumber\\
\left(\frac{\partial a_1}{\partial x^3}-\frac{\partial a_3}{\partial
x^1}\right)dx^3 dx^1$$
$$d\theta^2=\left(\frac{\partial a_{23}}{\partial
x^1}+\frac{\partial a_{31}}{\partial x^2}+\frac{\partial
a_{12}}{\partial x^3}\right)dx^1 dx^2 dx^3$$

One can see the following.

a) Any function is the form of zero degree. Its basis is a surface of zero
dimension, namely, a variety of points. The differential of the form of zero
degree is an ordinary differential of the function.

b) The form of first degree is a differential expression. As it was pointed
out above, an ordinary differential of a function is an example of the
first-degree form.

c) Coefficients of differentials of forms of zero-, first- and second degrees
give the gradient, curl, and divergence, respectively. That is, the operator
$d$, referred to as the exterior differentiation, is an abstract
generalization of the ordinary operators of gradient, curl,  and divergence.
At this point it should be emphasized that, in mathematical analysis the
ordinary concepts of gradient, curl, and divergence are the operators applied
to vectors, and in the theory of exterior forms, gradients, curls, and
divergences obtained as the results of exterior differentiating (forms of
the zero- first- and second degrees) are operators applied to 
pseudovectors (an axial vector).

These following examples are examples of exterior differential forms.

From these examples one can assure oneself that, firstly, the
differential of an exterior form is also an exterior form
(but with the degree greater by one), and, secondly, one can see
that the components of differential forms are commutators of the
coefficients of the form's differential. Thus, the differential
of the first-degree form $\omega=a_i dx^i$ can be written as
$d\omega=K_{ij}dx^i dx^j$ where $K_{ij}$ are the components of the
commutator for the form $\omega$ that are defined as
$K_{ij}=(\partial a_j/\partial x^i-\partial a_i/\partial x^j)$.

As pointed out, by definition an exterior differential form is
a skew-symmetric tensor field. To a differential form of
degree $p$ there corresponds a skew-symmetric covariant tensor
of the type $(0,p)$: $T=(T_{\alpha_1\dots\alpha_p})$. Such a
tensor may be written as
$$
T=\sum_{\alpha_1<\dots<\alpha_p}T_{\alpha_1\dots\alpha_p}e^{\alpha_1}\circ e^{\alpha_2}\circ\dots
\circ e^{\alpha_p}
$$
where $e^{\alpha_i}$ are base vectors. If the
differentials of coordinates $dx^{\alpha_i}$ are chosen as the basis,
then to the skew-symmetric tensor there will correspond the expression
$$
T=\sum_{\alpha_1<\dots<\alpha_p}T_{\alpha_1\dots\alpha_p}
dx^{\alpha_1}\wedge dx^{\alpha_2}\wedge\dots \wedge
dx^{\alpha_p}
$$
which is a differential form. (Since differentials of
coordinates must satisfy the condition of exterior multiplication (1.2),
the correspondence $dx^\alpha\wedge dx^\beta\leftrightarrow
e^\alpha e^\beta-e^\beta e^\alpha$ must be valid [3]. \{The
forms $dx^{\alpha_i}$ make up the basis of the cotangent space\}.

{\footnotesize [Historically, in physical applications the first method
proposed for expressing tensors, namely, by means of base vectors, was widely
used. Not only skew-symmetric tensors, but all tensors can be written in that
form. In the case when the
metric tensor exists, this way makes it possible to go from covariant indices
to contravariant ones, and vice verse. Nevertheless, a second method of
expressing the skew-symmetric tensors, i.e., as exterior differential forms,
has distinct advantages. (Practically all authors who work with exterior
differential forms direct attention to this fact [3, 5, 6, 7, 8]). The
method of presentation of the skew-symmetric tensors as differential forms
extends the capabilitites of the mathematical apparatus based on these
tensors. Tensors are known to have been introduced as objects that are
transform according to a fixed rule under transformation of coordinates.
Tensors are attached to a basis that can be transformed in an arbitrary way
under transition to a new coordinate map. An exterior differential form is
connected with differentials of coordinates that vary according to the
interior characteristics of the manifold, under translation along the
manifold. Using differentials of the coordinates, instead of the base
vectors, enables one to directly make use of the integration and
differentiation for physical applications. Instead of differentials
of the coordinates, a system of linearly-independent exterior one-forms can
be chosen as the basis, and this makes the description independent of the
choice of coordinate system [5, 6].]}

\subsection{Closed exterior differential forms}

In mathematical formalisms and mathematical physics 
the closed differential forms with
invariant properties appear to be of greatest practical utility. 

A form is called  `closed,'  if its differential is equal to zero:
$$
d\theta^p=0\eqno(1.4)
$$

From condition (1.4) one can see that a closed form is a conserved 
quantity. (This means that it corresponds to a conservation law, 
namely, to some conserved physical quantity.)

The differential of a form is a closed form. That is
$$
dd\omega=0\eqno(1.5)
$$
where $\omega$ is an arbitrary exterior form.

A form which is the differential of some other form: 
$$
\theta^p=d\theta^{p-1}\eqno(1.6)
$$
is called an `exact' form. Exact forms prove to be closed
automatically 
$$
d\theta^p=dd\theta^{p-1}=0\eqno(1.7)
$$
This follows from the property (1.5) of the the exterior differential
[5]. 

Here it is necessary to pay attention to the following points. In the
formulas presented above  it was implicitly assumed that the differential
operator $d$ is a total operator (i.e., it $d$ acts everywhere in the
vicinity of the point considered locally),  and therefore it acts on the
manifold of the initial dimension $n$. However, the differential may be
internal. Such a differential acts on some structure with the dimension being
less than that of the initial manifold. The structure, on which the exterior
differential form may become a closed inexact form, is a pseudostructure 
with respect to its metric properties. \{Cohomology
(de Rham cohomology, singular cohomology [4, 6]), sections of cotangent
bundles,  integra and potential surfaces and so on, may be regarded as
examples of pseudostructures. As it will be shown later, an eikonal surfaces
corresponds to a pseudostructure\}.  

If a form is closed on a pseudostructure only, the closure condition is
written as
$$
d_\pi\theta^p=0\eqno(1.8)
$$
(In this case the internal differential (rather then total one) becomes
equal to zero). And the pseudostructure $\pi$ obeys the condition
$$
d_\pi{}^*\theta^p=0\eqno(1.9)
$$
where ${}^*\theta^p$ is a dual form. (For the properties of dual forms, 
see [6]). 

{\footnotesize [In a two-dimensional space, $(x,y)$, an exterior differential 
form can be written as
$\theta =udx+vdy$, and the corresponding dual form is $^*\theta =-vdx+udy$.]}

From conditions (1.8) and (1.9) one can see that the form
closed on a pseudostructure is a conserved object, namely, this
quantity is conserved on a pseudostructure. (This can also correspond to
some conservation law, i.e. to a conserved object.)

An exact form is, by definition, a differential (see condition (1.6)).
In this case the differential is total. A
closed inexact form is a differential too; and, in this case the
differential is an interior one defined on a pseudostructure. Thus, any
closed form is a differential. The exact form is a total
differential. The closed inexact form is an interior (on
pseudostructure) differential, that is
$$
\theta^p_\pi=d_\pi\theta^{p-1}\eqno(1.10)
$$
At this point it is worth noting that the total differential of a
form closed on the pseudostructure is nonzero, that is
$$
dd_\pi\omega\ne0\eqno(1.11)
$$

And so, any closed form is a differential of a form of  lower
degree: the total one $\theta^p=d\theta^{p-1}$ if the form is exact,
or the interior one $\theta^p=d_\pi\theta^{p-1}$ on pseudostructure if
the form is inexact. (From this it follows that the form of lower
degree may correspond to a potential, and the closed form by itself
may correspond to a potential force. This is an additional example showing
that a closed form may have physical meaning. Here the two-fold nature
of a closed form is revealed, on the one hand, as a locally conserved
quantity, and on the other hand, as a potential force.)

From the conditions (1.6) and  (1.10), it is possible to see that
a connection between closed forms of different
degrees can exist. Due to  condition (1.5), for exact forms this
connection couples only two forms. If the form $\theta^p$ is
exact, then there exists a connection between the forms $\theta^p$
and $\theta^{p-1}$: $\theta^p=d\theta^{p-1}$, but the form
$\theta^{p+1}=d\theta^p$ vanishes according to condition
(1.11), and the connection is broken. In this case it is assumed that
the form $\theta^{p-1}$ is not exact, since otherwise the form
$\theta^p$ would be equal to zero. There is then no connection with the
form $\theta^{p-2}$. For closed inexact forms, for which
condition (1.11) is satisfied, such connections may couple greater
numbers of terms. \{In differential geometry this fact is
related to cohomology theory, the theory of structures\}.

Similarly to the differential connection between exterior
forms of sequential degrees, there is an integral connection. The relevant
integral relation has the form [6]
$$
\int\limits_{c^{p+1}}d\theta^p=\int\limits_{\partial
c^{p+1}}\theta^p\eqno(1.12)
$$
In particular, the integral theorems by Stokes and Gauss follow
from the integral relation for $p=1,2$ in three-dimensional space.
\{From this relation one can see that the integral of the closed form
over the closed curve vanishes (in the case of a smooth manifold).
However, in the case of a complex
manifold (for example, a not simply connected manifold, with the homology
class being nonzero), an integral of a closed form (in this case the
form is inexact) over a closed curve is nonzero. It may be equal to a
scalar multiplied by $2\pi$, which in this case corresponds to, for
example, a physical quantity such as charge [6]. Just such integrals
are considered in the theory of residues.\}

\subsection{Properties of the closed exterior forms} 

The role of closed exterior forms in matematics relates to the fact that 
the properties of closed exterior and dual forms, namely, invariance,
covariance, conjugacy, and duality, lie at the basis of the group,
structural and other invariant methods of matematics. 
 
\subsection*{Invariant properties of closed exterior differential
forms.} 

Since a closed form is a differential, then  it is
obvious that a closed form will turn out to be invariant under all 
transformations that conserve the differential. (The nondegenerate 
transformations in mathematics and mathematical physics such as the unitary, 
tangent, canonical, gradient, and other nondegenerate transformations 
are examples of such transformations that conserve the differential.) 

Invariant properties of closed exterior forms explicitly or implicitly
manifest themselves essentially in all invariant mathematical formalisms 
and formalisms of field theory, such as the Hamilton formalism, tensor 
calculus, group theory, quantum mechanics equations, Yang-Mills theory 
and others.

Covariance of a dual form is directly connected with the invariance of an
exterior closed form.

The invariance property of an closed inexact exterior form and covariance 
of its dual form play an important role in describing invariant structures 
and manifolds.

\bigskip
{\bf Invariant structures}

Of most significance in mathematical formalisms and mathematical physics 
are closed inexact exterior forms. This is due to the fact that 
the closed inexact exterior form and relevant dual form describe 
the differential-geometrical structure, which is invariant one. 

From the definition of a closed inexact exterior form one can see
that to this form there correspond two conditions:

(1) condition (1.8) is the closure condition of the exterior form 
itself, and

(2) condition (1.9) is that of the dual form.

Conditions (1.8) and (1.9) can be regarded as equations for a binary 
object that combines the pseudostructure (dual form) and 
the conserved quantity (the exterior differential form) defined
on this pseudostructure. Such a binary object is  
differential - geometrical structure.  (The well-known G-Structure
is an examlpe of such differential-geometrical structure.)  

As it has been already pointed out, closed inexact exterior form is  
a differential  (an interior one on the pseudostructure), and hence 
it remains invariant under all transforms that conserve the differential. 
Therefore, the relevant differential-geometrical structure also remains 
invariant under all transforms that conserve 
the differential. For the sake of convenience in the subsequent
presentation such differential - geometrical structures 
will be called the I-Structures.

To the unique role of such invariant structures in mathematics it points 
the fact that the transformations conserving the differential 
(unitary, tangent, canonical, gradient and gauge ones) lie at the basis 
of many branches of mathematics, mathematical physics and field theory. 

The invariant structures appear  
while analyzing the integrability of differential equations. 
Their role in the theory of differential equations relates to the fact 
that they correspond to generalized solutions which describe measurable 
physical quantities. In this case the integral surfaces with conservative 
quantities (like the characteristics, the characteristic surfaces, 
potential surfaces and so on) are invariant structures.
The examples of such studying the integrability of differential equations 
using the skew-symmetric differential forms are presented in paper [5].

The mechanism of realization of the differential-geometrical structures 
and their characteristics will be described in Subsection 2.4.

\subsection{Invariance as the result of conjugacy of elements of
exterior or dual forms}

Closure of exterior differential forms, and hence their invariance, results
from the conjugacy of the elements of exterior or dual forms. 

From the definition of an exterior differential form one can see that 
exterior differential forms have complex structure. The specific features of
the structure of exterior forma are homogeneity with respect to the basis,
skew-symmetry, the integration of terms each consisting of two objects of
different nature (the algebraic nature for the form coefficients, and the
geometric nature of the base components). Besides,  an exterior form depends
on the space dimension and on the manifold topology. The closure property of
an exterior form implies that any objects, namely, elements of the exterior
form, components of elements, elements of the form's differential, exterior
and dual forms and others, turn out to be conjugated. It is conjugacy
that leads to realization of the invariant and covariant 
properties of the exterior and dual forms that have great functional and
applied importance. The variety of conjugate objects leads to the fact that 
closed forms can describe a great number of different physical and spatial
structures, and this fact, once again, emphasizes the vast mathematical
capabilities of the exterior differential forms. 

Let us consider some types of conjugacy that make the exterior differential
and dual forms closed, that is, they make these form differentials equal to
zero.

As it was pointed out already, the components of an exterior form a
commutator are the
coefficients of the differential of this form. If the commutator of the form
vanishes, the form differential vanishes too, and this indicates that the
form is closed. Therefore, closure of a form may be recognized by finding
whether or not the commutator of the form vanishes.

One of the types of conjugacy is that for the form coefficients.

Let us consider an exterior differential form of the first degree
$\omega=a_i dx^i$. In this case the differential will be expressed
as $d\omega=K_{ij}dx^i dx^j$, where
$K_{ij}=(\partial a_j/\partial x^i-\partial a_i/\partial x^j)$ are
the components of the form's commutator.

It is evident that the differential may vanish if the components
of commutator vanish. One can see that
the components of the commutator $K_{ij}$ may vanish if derivatives of
the form's coefficients vanish. This is a trivial
case. In addition, the components $K_{ij}$ may vanish if the
coefficients $a_i$ are derivatives of some function $f(x^i)$,
that is, $a_i=\partial f/\partial x^i$. In this case,
the components of the commutator are equal to the difference of the mixed
derivatives
$$
K_{ij}=\left(\frac{\partial^2 f}{\partial x^j\partial
x^i}-\frac{\partial^2 f}{\partial x^i\partial x^j}\right)
$$
and therefore they vanish. One can see that those form
coefficients $a_i$, that satisfy these conditions,  are
conjugated quantities (the operators of mixed differentiation turn out
to be commutative).

Let us consider the case when the exterior form is written as
$$
\theta=\frac{\partial f}{\partial x}dx+\frac{\partial f}{\partial
y}dy
$$
where $f$ is the function of two variables $(x,y)$. It is evident that this
form is closed because it is
equal to the differential $df$. And for its dual form
$$
{}^*\theta=-\frac{\partial f}{\partial y}dx+\frac{\partial
f}{\partial x}dy
$$
to be closed also,  it is necessary that its commutator be equal to zero
$$
\frac{\partial^2 f}{\partial x^2}+\frac{\partial^2 f}{\partial
y^2}\equiv \Delta f=0
$$
where $\Delta$ is the Laplace operator. As a result, the function $f$ has
to be harmonic.

Let us assume that the exterior differential form of first degree has
the form $\theta=udx+vdy$, where  $u$ and $v$ are functions of
two variables $(x,y)$. In this case, the closure condition of the
form, that is, the condition under which the form commutator vanishes,
takes the form
$$
K=\left(\frac{\partial v}{\partial x}-\frac{\partial u}{\partial
y}\right)=0
$$
One can see that this is one of the Cauchy-Riemann conditions for complex
functions. The closure condition for the relevant dual form
${}^*\theta=-vdx+udy$ is the second Cauchy-Riemann condition. \{Here one can
see a connection between exterior differential forms and  functions of
complex variables. If we consider the function $w=u+iv$ of the complex
variables $z=x+iy$, $\overline{z}=x-iy$ which satisfy  the Cauchy-Riemann
conditions, then to closed exterior and dual forms will correspond to this
function. (The Cauchy-Riemann conditions are the conditions under which a
function of complex variables does not depend on the conjugated coordinate
$\overline{z}$). And, to each harmonic function of complex variables there
corresponds the closed exterior differential form, whose coefficients $u$ and
$v$ are conjugated harmonic functions\}.

A conjugacy, which makes an interior differential on apseudostructure equal
to zero, can exist, say, $d_\pi\theta=0$. Let us assume that an interior
differential is a form of first degree (the form itself is that of zero
degree), and it can be given in the form $d_\pi\theta=p_x dx+p_y dy=0$, where
$p$ is a form of zero degree (any function). Then the closure condition of
the form is
$$
\frac{dx}{dy}=-\frac{p_y}{p_x}\eqno(1.13)
$$
This is a conjugacy of the basis and the derivatives of the form's coefficients. 
One can see that formula (1.13) is one of the formulas of the canonical
relations [11]. The second formula of the canonical relations follows from
the condition that the dual form differential vanishes. This type of
conjugacy  is connected with a canonical transformation. For the differential of 
first degree form (in this case the differential is a form of second
degree), the corresponding transformation is a gradient one. 

Relation (1.13) is the condition for an implicit function to exist. 
That is, the closed (inexact) form of zero degree is the implicit function.
\{This is an example  of a connection between an exterior forms and
analysis\}.

The property of the exterior differential forms being closed on a
pseudostructure points to another type of conjugacy, namely, the conjugacy of
exterior and dual forms.

\subsection*{Duality of the exterior differential
forms}
Conjugacy is connected with another characteristic property of exterior
differential forms, namely, their duality; and, as will be shown later, it
has fundamental physical meaning. \{Conjugacy is an identical connection
between two operators or mathematical objects.  Duality is a concept meaning
that one object carries a double meaning, or that two objects with different
meanings (of different physical nature) are identically connected. If one
knows any dual object, one can obtain the other object\}. The conjugacy of
objects of the exterior differential form generates duality of the exterior
forms. 

The connection between exterior and dual forms is an example of  duality. An
exterior form and its dual form correspond to  objects of different nature:
an exterior form corresponds to a physical (i.e. algebraic) quantity, and its
dual form corresponds to  some spatial (or pseudospatial) structure. At the
same time, under conjugacy, the duality of these objects manifests itself,
that is, if one form is known it is possible to find the other form. \{It
will be shown below that the duality between exterior and dual  forms
elucidates a connection between physical quantities and spatial  structures
(together they form a physical structure\}. Here duality is also evident in
the fact that, if the degree of  the exterior form equals $p$, the dimension
of the structure equals $N-p$, where $N$ is the space dimension.

Since a closed exterior form possesses invariant properties and the dual form
corresponding to it possesses covariant properties,  invariance of the closed
exterior form and the relevant covariance of its dual form, is an example of
the duality of exterior differential forms.

Another example of duality of closed forms is connected with the fact that a
closed form of degree $p$ is a differential form of degree $p-1$. This
duality is manifested in that, on the one hand, as it was pointed out before,
the closed exterior form is a conserved quantity, and on the other hand, the
closed form can correspond to a potential force. (Below the physical meaning
of this duality will be elucidated, and it will be shown in respect to what a
closed form manifests itself as a potential force and with what the
conserved physical quantity is connected).

A further manifestation of  duality of exterior differential forms, that needs
more attention, is the duality of the concepts of closure and of
integrability of differential forms.

A form that can be presented as a differential can be called integrable,
because it is possible to integrate it directly. In this context a closed
form, which is a differential (total, if the form is exact, or interior, if
the form is inexact), is integrable. And the closed exact form is integrable
identically, whereas the closed inexact form is integrable on pseudostructure
only.

The concepts of closure and integrability are not identical. Closure of a
form is defined with respect to the form of degree greater by one, whereas
integrability is defined with respect to a form of degree less by one.
Really, a form is closed if the form's differential, which is a form of
greater by one degree, equals zero. And a form referred to as integrable is a
differential of some form of degree less by one. Closure and  integrability
are dual concepts. Namely, closure and  integrability are further examples of
duality among exterior differential forms.

Here it should be emphasized that duality is a property of closed
forms only. In particular, nonclosure and nonintegrability are not dual
concepts, and this will be revealed while analyzing the evolutionary
differential forms.

Duality of closed differential forms is revealed by an availability of one or
another type of conjugacy. As it will be shown below, this has physical
significance. Duality is a tool that untangles the mutual connection, the
mutual changeability and the transitions between different physical objects
in the evolutionary processes.

\subsection{Specific features of the mathematical apparatus of exterior
differential forms}

\subsection*{Operators of the theory of exterior differential
forms}

The distinguishing properties of the mathematical apparatus of exterior
differential forms were identified by Cartan [7]. 
Cartan aimed to build a theory, which contains concepts
and operations being independent of any change of variables, both
dependent and independent. To do this it was necessary to exchange partial
derivatives for differentials that have an interior meaning.

In  the differential calculus, derivatives are the basic elements of the
mathematical apparatus. By contrast, a  differential is an element of the
mathematical apparatus of the theory of exterior differential forms. It
enables one to analyze the conjugacy of derivatives in various
directions, which extends the potentialities of the differential
calculus. (A derivative that can be considered as the conjugacy of the
differentiating operator and scalar function is an analog of the exterior
zero-degree differential).

The exterior differential operator, $d$, is an abstract generalization  of
the ordinary mathematical operations of gradient, curl and divergence 
from vector calculus [5]. If, in addition to the exterior differential, 
we introduce the following operators: (1) $\delta$ for transformations that 
convert the form of degree $p+1$ into the form of degree $p$, (2) 
$\delta'$ for cotangent transformation, (3) $\Delta$ for the $d\delta-\delta
d$ transformation, (4)$\Delta'$ for the $d\delta'-\delta'd$ transformation, then in 
terms of these operators, that act on the exterior differential forms,
one can write down the operators in field theory equations. The operator
$\delta$ corresponds to Green's operator, $\delta'$ to the canonical
transformation operator, $\Delta$ to the d'Alembert operator in 4-dimensional
space, and $\Delta'$ corresponds to the Laplace operator [6, 8]. It
can be seen that the operators of the exterior differential form theory
are connected with many operators of mathematical physics.

The mathematical apparatus of exterior differential forms extends the
capabilitites of integral calculus. As pointed out, exterior differential
forms were introduced as integrand expressions for definition of the integral 
invariants. The closure condition for exterior differential forms makes
it possible to find  integrability conditions, namely, to find the
conditions under which the integral is independent of the integration
curve. The connection between forms of subsequent degrees establishes the
relations which under some conditions make it possible to reduce
integration over some domain to that along some boundary, and vice verse
(see formula (1.6)). This fact is of great importance for applications.

\subsection*{Identical relations of exterior differential forms} 

In the theory of exterior differential forms closed forms, that possess
various types of conjugacy, play a principal role. Since conjugacy
represents a certain connection between two operators or mathematical
objects, it is evident that relations can  be used to express conjugacy
mathematically. Just such relations constitute the basis of the
mathematical apparatus of exterior differential forms.

At this point the following it should be emphasized. A relation is a
comparison, a correlation of two objects. A relation may be an
identical or nonidentical (see Section 2). 

The basis of the mathematical apparatus of exterior differential forms
comprises identical relations. (Below nonidencal relations
will be discussed, and it will be shown that identical relations for
exterior differential forms are obtained from nonidencal relations.
Also it will be shown, that transitions from nonidentical relations to
identical ones describe transitions of any quality into another one).

The identical relations of exterior differential forms reflect the
closure conditions of the differential forms, namely, the vanishing of
the form's differential (see formulas (1.4), (1.8), (1.9)) and the
conditions connecting forms of subquential degrees (see formulas (1.6),
(1.10), (1.12)).

The importance of identical relations for exterior differential forms
is manifested by the fact that practically in all branches of physics,
mechanics, thermodynamics one encounters such identical relations.

Several kinds of identical relations can be distinguished.

1. \emph{ Relations among differential forms}.

These are  relations connecting forms of sequental degrees. They
correspond to formulas (1.6), (1.10). Examples of such identical
relations are

a) the Poincare invariant $ds\,=\,-H\,dt\,+\,p_j\,dq_j$,  

b) the second principle of thermodynamics $dS\,=\,(dE+p\,dV)/T$,

c) the vital force theorem in theoretical mechanics: $dT=X_idx^i$
where $X_i$ are the components of a potential force, and $T=mV^2/2$ is the
vital force,

d) the conditions on characteristics [9] in the theory of differential
equations.

The requirement that the function is an antiderivative (the integrand is a
differential of a certain function) can be written in terms of such an 
identical relation.

The existence of a harmonic function is expressed by means  of an
identical relation: a harmonic function is a closed form, that is, a
differential (a differential on the Riemann surface).

Identical relations among differential forms expresses the fact that
each closed exterior form is a differential of some exterior form 
(with the degree less by one). 

In general, such an identical relation can be written as
$$
d\phi=\eta^p\eqno(1.14)
$$

In this relation the form on the right-hand side has to be a \emph{
closed.} ( As will be shown below, the identical relations are satisfied only on
pseudostructures. That is,  an identical relation can be written as  
$$ 
d_{\pi}\phi=\eta_{\pi}^p\eqno(1.14') 
$$

In the identical relations (1.14) on one side there is a closed form, and
on other  a differential of some differential form of degree the less by
one.

From such identical relations one can obtain the following information.

If there is a differential of some exterior form, it means that there is
a closed exterior form (which points, for example, to the availability
of the structure).

On the other hand, the availability of a closed exterior form points to 
the availability of a differential (and this may mean that there is 
a potential or a state function). 

From such a relation it follows that a differential is the result of
conjugacy of some objects, because this differential is a closed form.
Hence a potential or a state function is also connected with the
conjugacy of certain objects.

2. \emph{ Integral identical relations}.

At the beginning of the paper, it was pointed out that exterior
differential forms were introduced as integrands possessing the
following property: they can have integral invariants. This fact (the
availability of an integral invariant) is mathematically expressed  as
an identical relation.

Many formulas of Newton, Leibnitz, Green, the integral relations by 
Stokes, Gauss-Ostrogradskii are examples of integral identical relations 
(see formula (1.12)).

3. \emph{ Tensor identical relations}.

From relations that connect exterior forms of sequential degrees, one
can obtain vector and tensor identical relations that connect the
operators: gradient, curl, divergence and so on.

From closure conditions on exterior and dual forms, one can obtain 
identical relations such as the gauge relations in electromagnetic field
theory, the tensor relations between connectednesses and its derivatives
in gravity theory (the symmetry of connectednesses with respect  to
lower indices, the Bianchi identical, the conditions imposed on the 
Christoffel symbols), and so on.

4. \emph{ Identical relations between derivatives}.

Identical relations between derivatives correspond to closure
conditions on exterior and dual forms. Examples of such relations 
are the above presented Cauchi-Riemann conditions in the theory of 
complex variables, the transversality condition in the calculus of 
variations, the canonical relations in the Hamilton formalism, the 
thermodynamic relations between derivatives of thermodynamic functions 
[10], the condition that the derivative of implicit functions are
subject to, the eikonal relations [11], and so on.

The examples presented above show that identical relations between
exterior differential forms occur in various branches of mathematics and
physics.

Application of identical relations, provided with a knowledge of closed
forms, provides the capability (a) to find other forms that are
necessary for describing physical phenomena, (b) to answer the question
of whether the obtained exterior forms are closed, (c) to get
information conveyed by the closed forms, and so on. They make it
possible to find closed differential forms that  correspond to
conservation laws, invariant structures, structures of manifolds and so on. 
They allow potentials and state functions to be determined.

Identical relations between exterior differential forms are  mathematical
expressions of various kinds of conjugacy that lead to closed exterior
forms. They describe the conjugacy of many objects: the form elements,
components of each element, exterior and dual forms, exterior forms of
various degrees, and others. Identical relations that are connected with
different kinds of conjugacy (the form elements, components of each
element, forms of various degrees, exterior and dual forms, and others)
elucidate invariant, structural and group properties of exterior forms,
which are of great importance in applications. The invariant, structural
and group properties of exterior forms that have physical meaning,
manifest themselves  as a result of one or another kind of conjugacy.
This is expressed mathematically as an identical relation.

The functional significance of identical relations for
exterior differential forms lies in the fact that they can describe the
conjugacy of objects of different mathematical nature. This enables one
to see the internal connections between various branches of mathematics.
Due to these possibilities, exterior differential forms have wide
application in various branches of mathematics. (Below a connection of
exterior differential forms with different branches of mathematics will
be discussed).

The mathematical apparatus of exterior differential forms is based on
identical relations of great utilitarian and functional importance.
However, the question of how closed exterior differential forms and 
identical relations are
obtained, and how the process of conjugating different objects is
realized, remains open.

The mathematical apparatus of evolutionary differential forms that is based
on nonidentical relations enables us to answer this question.
This will be shown in Section 2.

\subsection*{Nondegenerate transformations}

One of the fundamental methods in the theory of exterior differential
forms is application of \emph{ nondegenerate} transformations 
(below it will be shown that \emph{ degenerate} transformations appear 
in the mathematical apparatus of the evolutionary forms).

In the theory of exterior differential forms nondegenerate 
transformations are those that conserve the differential. 
Unitary transformations 
(0-forms), the tangent and canonical transformations (1-forms), gradient and 
gauge transformations (2-forms) and so on are  examples of such 
transformations.  These are the gauge transformations
for spinor, scalar, vector, tensor  (3-form) fields.
                                                       
From the descriptions of operators comprised of exterior differential
forms, one can see that they are operators that execute some
transformation. All these transformations are connected with the above listed
nondegenerate transformations of exterior differential forms.
(It is worth pointing out that just such transformations are used extensively 
in field theory. The field theory operators such as Green's, d'Alembert's and
Laplace operators are connected with nondegenerate transformations of closed 
exterior differential forms.)

The possibility of applying nondegenerate transformations shows that 
exterior differential forms possess group properties. This extends the
utilitarian potentialities of the exterior differential forms.

The significance of nondegenerate transformations consists in the fact, 
that they allow one to get new closed differential forms, which
opens up the possiblity of obtaing new structures.

Nondegenerate transformations, if applied to identical relations, enable one
to obtain new identical relations and new closed exterior differential 
forms.

Specific features of nondegenerate transformations, and their relation to the
degenerate transformations, will be discussed in Subsection 2.3 in more detail.

\bigskip

From the properties of exterior differential forms presented above one
can see, that the properties of exterior differential forms conform to
specific features of many branches of mathematics. Below we outline
connections between exterior differential forms and various branches of
mathematics in order to show what role exterior differential forms play
in mathematics, and to draw attention to their great potentialities.

\subsection{Connection between exterior differential forms
and various branches of mathematics}

\subsection*{Connection with tensors}

In Section 1.1 a connection between exterior differential forms and 
skew-symmetrical tensors was demonstrated. More detailed information on
this subject can be found, e.g., in [3,5,8].

As was pointed out in Section 1.1, the exterior differentiation
operator, $d$, is the abstract generalization of the gradient, curl and
divergence. This property elucidates a connection between vectorial,
algebraic and potential fields. The property of the exterior
differential form, namely, the existence of differential and integral
relations between the forms of sequential degrees, allows us to classify
these fields according to their exterior  form degree. 

\subsection*{Algebraic properties of exterior differential
forms}

The basis of Cartan's method of exterior forms, namely, the method of
analyzing systems of differential equations and manifolds, is the basis
of the `Grassmann' algebra (i.e., exterior algebra) [1]. The
mathematical apparatus of exterior differential forms extends  algebraic
mathematical techniques. Differential forms treated as elements of
algebra, allow studying manifold structure and finding  manifold
invariants. The group properties of exterior differential forms (the
connection with the Lie groups) enable one to study the integrability of
differential equations. They can constitute the basis of the invariant
field theory. They establish a connection between vectorial and
algebraic fields.

\subsection*{Geometrical properties of exterior differential forms}

Exterior differential forms elucidate the internal connection between
algebra and geometry. From the definition of a closed inexact form it
follows that a closed inexact form is a quantity that is conserved on a
pseudostructure. That is, a closed inexact form is a conjugacy of the
algebraic and geometrical approaches. The set of relations, namely, the
closure condition for exterior forms, ($d_\pi \theta^p=0$), and the
closure condition for dual forms, ($d_\pi {}^*\theta^p=0$), allow the
description of conjugacy. Closed form possess algebraic (invariant)
properties, and  closed dual forms have geometrical (covariant)
properties.

The mathematical apparatus of exterior differential forms enables one to
study the elements of the interior geometry. This is a fundamental
formalism in differential geometry [7]. It enables one to investigate
manifold structure. Cartan developed the method of exterior forms for
the purpose of investigating manifolds. It is known that in the case of
integrable manifolds, the metrical and differential characteristics are
consistent. Such manifolds may be regarded as objects of the interior
geometry, i.e. the possibility of studying their characteristics as the
properties of the surface itself without regard to its embedding space.
Cartan's structure equation is a key tool in studying manifold
structures and fiber spaces [7,12].

\subsection*{Theory of functions of complex variables}

The residue method in the theory of analytical functions of complex
variables is based on the integral theorems by Stokes and
Cauchy-Poincare that allow us to replace the integral of closed form
along any closed loop by the integral of this form along another closed
loop that is homological to the first one [2].

As already noted, harmonic functions (differentials on Riemann surfaces)
are closed exterior forms:
$$
\theta=pdx+qdy,\quad d\theta=0
$$
to which there corresponds the dual form:
$$
{}^*\theta=-qdx+pdy,\quad d{}^*\theta=0
$$
where $d\theta=0$ and $d{}^*\theta=0$ are the Cauchy-Riemann
relations.

\subsection*{Differential equations}

On the basis of the theory of exterior differential forms the methods
for studying the integrability of a system of differential equations,
Pfaff equations (the Frobenius theorem), and of finding integral
surfaces have been  developed. This problem was considered in many works
concerning exterior differential forms [5, 7, 13].

An example of an application of the theory of exterior differential
forms for analyzing integrability of differential equations and
determining  the functional properties of the solutions to these
equations is  presented in Subsection 1.6.

\subsection*{Connection with the differential and integral calculus}

As it was already pointed out, exterior differential forms were
introduced for designation of integrand expressions which can form
integral invariants [1].

Exterior differential forms are connected with multiple integrals (see,
for example [3]).

It should be mentioned that the closure property of the form $f(x)dx$
indicates existence of an antiderivative of the function $f(x)$.

Integral relations (1.12), from which the formulas by Newton, Liebnitz,
Green, Gauss-Ostrogradskii, Stokes were derived, are of great
importance. Such potentials as Newton's, the  potentials of simple and
double layers are integrals of closed inexact forms.

Exterior differential forms extend the potentialities of the
differential and integral calculus.

The theory of integral calculus establishes a connection between
the differential form calculus and the homology of manifolds
(cohomology).

The operator $d$ appears to be useful for expressing the integrability
conditions for systems of partial differential equations.

\bigskip 

Thus, even from this brief description of the properties and specific
features of the exterior differential forms and their mathematical
apparatus, one can clearly see their connection with such branches of
mathematics as  algebra, geometry, mathematical analysis, tensor
analysis, differential geometry, differential equations, group theory,
theory of transformations and so on. This is indicative of wide functional
and utilitarian potentialities of exterior differential forms.  Exterior
differential forms enable us to see the internal connection  between
various branches of mathematics and physics.

Here it arises the question of how closed exterior differential  forms,
which possess vast potentialities, are obtained. 

It appears, that closed exterior differential forms can be obtained 
from the evolutionary differential forms (see Section 2). 
 
Below, an example will be presented of an application of the
mathematical apparatus of exterior differential forms that demonstrates
the existence of  new evolutionary skew-symmetric differential forms. 

\subsection{Qualitative investigation of the functional properties
of the solutions to differential equations}

The investigation of the functional properties of the solutions to
differential equations on the basis of the mathematical apparatus of
skew-symmetric differential forms  makes it possible to understand what 
lies at the basis of the qualitative theory of differential equations.
(The presented method of investigating solutions to differential 
equations is not new. Such an approach was developed already by Cartan
[7] in his analysis of the integrability of differential equations. Here
it is featured to  demonstrate the significance of new skew-symmetric
differential forms.)

The role of exterior differential forms in the qualitative investigation 
of solutions to differential equations is a consequence of the fact  that
the mathematical apparatus of these forms enables one to determine  the
conditions for consistency of various elements of differential equations
or of systems of differential equations. This enables one, for example, to
define the consistency of  partial derivatives in  partial differential
equations, the consistency of  differential  equations in  system of
differential equations, the conjugacy of the  function derivatives and of
derivatives of initial data for ordinary  differential equations and so
on. Functional properties of  solutions to differential equations depend
just on whether or not  conjugacy conditions are satisfied.

The basic idea of the qualitative investigation of the solutions to
differential equations can be clarified by the example of the first-order
partial differential equation.

Let
$$ F(x^i,\,u,\,p_i)=0,\quad p_i\,=\,\partial u/\partial x^i \eqno(1.15)$$
be the partial differential equation of the first order. 

Let us consider the functional relation
$$ du\,=\,\theta\eqno(1.16)$$
where $\theta\,=\,p_i\,dx^i$ (the summation over repeated indices is 
implied). Here $\theta\,=\,p_i\,dx^i$ is a differential form of the first
degree.

The specific feature of a functional relation (1.16) is that in
the general case this relation turns out to be a nonidentical. 

The left-hand side of this relation involves a differential, and the
right-hand side includes a differential form  $\theta\,=\,p_i\,dx^i$. For
this  relation to be an identical, the differential  form
$\theta\,=\,p_i\,dx^i$ must be a differential as well (like the left-hand
side of relation (1.16)), that is, it has to be a closed exterior
differential form. To achieve this requires the commutator
$K_{ij}=\partial p_j/\partial x^i-\partial p_i/\partial x^j$ of the
differential form $\theta $, vanish.

In the general case, from equation (1.15) it does not follow (explicitly) 
that the derivatives $p_i\,=\,\partial u/\partial x^i $ satisfy the
equation (and given boundary or initial conditions of the problem) and
constitute a differential. Without supplementary conditions, the
commutator of the differential form $\theta $ defined as $K_{ij}=\partial
p_j/\partial x^i-\partial p_i/\partial x^j$ is not equal  to zero. The
form $\theta\,=\,p_i\,dx^i$ proves to be unclosed and is not a a
differential like the left-hand side of relation (1.16). The functional
relation (1.16) appears to be a nonidentical: the left-hand side of this
relation is a differential, but the right-hand side is not a differential.

The nonidentity of the functional relation (1.16) points to the fact, that
without additional conditions the derivatives of the initial equation do
not constitute a differential. This means that the corresponding solution
to the differential equation $u$ will not be a function of $x^i$. It will
depend on the commutator of the form $\theta $, that is, it will be a
functional.

To obtain the solution that is the function (i.e., the derivatives of this
solution made up a differential), it is necessary to add a closure
condition for the form $\theta\,=\,p_idx^i$ and for its dual form (in the
present case the functional $F$ plays a role of the form dual to 
$\theta$) [7]:
$$\cases {dF(x^i,\,u,\,p_i)\,=\,0\cr
d(p_i\,dx^i)\,=\,0\cr}\eqno(1.17)$$
If we expand the differentials, we get a set of homogeneous equations with
respect to $dx^i$ and $dp_i$ (in the $2n$-dimensional space -- initial and
tangential):
$$\cases {\displaystyle \left ({{\partial F}\over {\partial x^i}}\,+\,
{{\partial F}\over {\partial u}}\,p_i\right )\,dx^i\,+\,
{{\partial F}\over {\partial p_i}}\,dp_i \,=\,0\cr
dp_i\,dx^i\,-\,dx^i\,dp_i\,=\,0\cr} \eqno(1.18)$$
The solvability conditions for this set (i.e., vanishing of the
determinant composed of the coefficients at $dx^i$, $dp_i$) have the form:
$$
{{dx^i}\over {\partial F/\partial p_i}}\,=\,{{-dp_i}\over 
{\partial F/\partial x^i+p_i\partial F/\partial u}} \eqno (1.19)
$$
These conditions determine an integrating direction, namely, a
pseudostructure,  on which the form $\theta \,=\,p_i\,dx^i$ turns out to
be closed, i.e. it becomes a differential, so that from relation (1.16)
the identical relation is deduced. If conditions (1.19), that may be
called `integrability conditions,' are satisfied, the derivatives
constitute a differential  $\delta u\,=\,p_idx^i\,=\,du$ (on the
pseudostructure), and the solution becomes a function. Just such
solutions, namely, functions on the pseudostructures formed by the
integrating directions, are the so-called generalized solutions [14].
Derivatives of generalized solutions made up the exterior forms that are
closed on pseudostructures.

(In this case the pseudostructure, i.e. the integral curve, and the 
dufferential $du$, which is a closed form, made up an I-Structure. 
As one can see, the realization of I-Structure, i.e. 
the realization of conditions (1.19) and closed form, are necessary 
conditions for constructions the generalized solution. 
The generalized solution is a solution of the
equation that correaponds  to the moment of realization of the conditions
of degenerate transformation.  The pseudostructure, on which such a
minisolution is defined, is  correspondingly a minisrtucture. Such
minisolutions are realized in discrete manner.)  

If conditions (1.19) are not satisfied, that is, the derivatives do not
form a differential, then the solution that corresponds to such
derivatives will depend on the differential form commutator consisting of
derivatives. That means that the solution is a functional rather then a
function. 

Since functions that are the generalized solutions are defined only on
pseudostructures, they have discontinuities in their derivatives in
directions that are transverse to the pseudostructures. The order of the 
derivatives with discontinuities equals the  degree of the exterior form.
If a form of zero degree is involved in the functional relation, then
the function itself, being a generalized solution, will have
discontinuities.

If we find the characteristics of equation (1.15), it appears that
conditions (1.19) are the equations for the characteristics [9]. That is,
the characteristics are examples of the pseudostructures on which the
derivatives of the differential equation consist of closed forms and the
solutions prove to be the functions (generalized solutions). (The
characteristic manifolds of equation (1.15) are the pseudostructures $\pi$
on which the form $\theta =p_idx^i$ becomes a closed form: $\theta
_{\pi}=d u_{\pi}$).
In this case the characteristics are pseudostructures. The characteristics and 
relevant closed forms on characteristics made up I-Structures.

Here it is worth noting, that coordinates of the equations of
characteristics are not identical to the independent coordinates of the
initial manifold on which equation (1.15) is defined. The transition from the
initial manifold to the characteristic manifold appears to be a \emph{
degenerate} transform, namely, the determinant of the set of equations
(1.18) becomes zero. The derivatives of equation (1.15) are transformed
from the tangent manifold to the cotangent manifold. The transition from the
tangent manifold, where the commutator of the form $\theta$ is nonzero (the
form is not closed, and derivatives do not form a differential), to the
characteristic manifold, namely, the cotangent manifold, where the commutator
becomes equal to zero (and a closed exterior form arises, i.e. the
derivatives form a differential), is an example of a degenerate transform.

A partial differential equation of the first order has been analyzed, and
the functional relation with the form of the first degree analogous to an
evolutionary form has been considered.

Similar functional properties are possesed by solutions to all
differential  equations describing any processes. And, if the order of the
differential equation is $k$, a functional relation comprising  $k$-degree
form corresponds to this equation. For ordinary differential equations the
commutator is generated at the expense of the conjugacy of the derivatives
of the functions desired and those of the initial data (the dependence of
the solution on the initial data is governed by the commutator).

In a similar manner, one can also investigate  solutions to a set of
partial differential equations. 

It can be shown that  the solutions to equations of mathematical physics,
on which no additional external conditions are imposed, are functionals.
The solutions prove to be exact (the generalized solution)  only under
realization of additional requirements, namely,  conditions om degenerate
transforms: vanishing the determinants, Jacobians and so on, that define
the integral surfaces.  The characteristic manifolds, the envelopes of
characteristics, singular points, potentials of simple and double layers,
residues and others are the examples of such surfaces.

Here mention should be made of the generalized Cauchy problem when 
initial conditions are given on some surface. The so called ``unique"
solution to the Cauchy problem, when the output derivatives can be
determined (that is, when the determinant built of the expressions at
these derivatives is nonzero), is a functional since the derivatives
obtained in this a way turn out to be nonconjugated, that is, their mixed
derivatives form a nonzero commutator, and the solution depends on this
commutator.

The dependence of the solution on the commutator may lead to instability
of the solution. Equations that do not satisfy  integrability conditions
(the conditions such as, for example, characteristics, singular points,
integrating factors and others) may have unstable solutions. Unstable
solutions appear in the case when the additional conditions are not
realized and no exact solutions (their derivatives form a differential)
are formed. Thus, the solutions to the equations of the elliptic type may
be unstable. 

Investigation of nonidentical functional relations lies at the basis of the
qualitative theory of differential equations. It is well known that the
qualitative theory of differential equations is based on the analysis of
unstable solutions and the integrability conditions. From the functional
relation it follows that the dependence of the solution on the commutator
leads to instability, and the closure conditions of skew-symmetric forms
constructed by derivatives are integrability conditions. That is,
the qualitative theory of differential equations that solves the problem of
unstable solutions and integrability bases on the properties of nonidentical
functional relation.

Thus, application of the exterior differential forms allows one to expose
the functional properties of the solutions to differential equations.

Here it is possible to see two specific features that lie beyond the scope
of the exterior differential form theory: the appearance of a nonidentical
relation (see, functional relation (1.16)) and a degenerate transformation. 
The skew-symmetric differential forms, which, in contrast to exterior 
forms, are defined on nonintegrable manifolds (such, for example,
as the tangent manifolds of differential equations) and possess
the evolutionary properties, are provided with such mathematical
apparatus.

\section{Evolutionary skew-symmetric differential forms}

In this Section the evolutionary skew-symmetric  differential
forms are described. Such skew-symmetric forms are deduced from 
evolutionary differential equations and, in contrast to exterior forms, 
are defined on nonintegrable manifolds (such as tangent manifolds of 
differential equations, Lagrangian manifolds and so on). 
A specific feature of evolutionary forms is the fact that from 
the evolutionary forms the closed exterior forms are obtained.

The process of extracting closed exterior forms from evolutionary
forms describes the generation of various invariant structures, conjugated 
objects and operators, discrete transitions and quantum jumps, forming 
pseudometric and metric manifolds.

A distinction of evolutionary skew-symmetric differential forms from
exterior forms is connected with the properties of manifolds on which
skew-symmetric forms are defined.

In the beginning of this section information on skew-symmetric 
differential forms and the manifolds on which skew-symmetric differential 
forms can be defined is presented.

\subsection{Some properties of manifolds}

Assume that on a manifold one can place a coordinate system with base
vectors $\mathbf{e}_\mu$ and define the metric forms for a manifold [15]:
$(\mathbf{e}_\mu\mathbf{e}_\nu)$, $(\mathbf{e}_\mu dx^\mu)$,
$(d\mathbf{e}_\mu)$.

If a metric form is closed (i.e., its commutators equal zero), then this
metric is defined by $g_{\mu\nu}=(\mathbf{e}_\mu\mathbf{e}_\nu)$ and the
results of a translation over a manifold of the point
$d\mathbf{M}=(\mathbf{e}_\mu dx^\mu)$ and of the unit frame
$d\mathbf{A}=(d\mathbf{e}_\mu)$ prove to be independent of the the path of
integration. Such a manifold is integrable. [On the integrability of
manifolds, see[15]].

If metric forms are nonclosed (the commutators of metric forms  are
nonzero), this points to the fact that this manifold  is nonintegrable.

Metric forms and their commutators define the metric and  differential
characteristics of a manifold.

Closed metric forms define a manifold structure, i.e. the internal 
characteristics of a manifold. And, nonclosed metric forms define  the
differential characteristics of a manifold. The topological properties 
of manifolds are connected with commutators of nonclosed metric forms. 
The commutators of nonclosed metric forms define the manifold 
differential characteristics that specify the manifold deformations: 
bending, torsion, rotation, twist.  Thus, the final result is, that
nonintegrable manifolds, i.e. the manifolds with  nonclosed metric forms, 
are deformed manifolds.

To describe manifold differential characteristics and, correspondingly, 
metric form commutators, one can use connectedness [2,5,7,15].

If the components of a metric form can be expressed in terms of
connectedness $\Gamma^\rho_{\mu\nu}$ [15], the expressions
$\Gamma^\rho_{\mu\nu}$, $(\Gamma^\rho_{\mu\nu}-\Gamma^\rho_{\nu\mu})$ and
$R^\mu_{\nu\rho\sigma}$ are components of the commutators of metric forms
of zeroth- first- and third degrees.  (The commutator of the second
degree metric form is written down in  a more complex manner [15], and
therefore it is not presented here).

As is known, a commutator of the zeroth degree metric form
$\Gamma^\rho_{\mu\nu}$ characterizes the bend, while that of the
first degree form $(\Gamma^\rho_{\mu\nu}-\Gamma^\rho_{\nu\mu})$
characterizes the torsion, the commutator of the third degree metric form
$R^\mu_{\nu\rho\sigma}$ determines the curvature. 

In the case of nonintegrable manifolds, 
the components of the metric form
commutators are nonzero.   In particular, the connectednesses 
$\Gamma^\rho_{\mu\nu}$ are not  symmetric. (For manifolds with a closed
metric form of the first degree,  the connectednesses  are
symmetric.)

Examples of nonintegrable manifolds are 
the tangent manifolds of differential equations that describe 
an arbitrary processes, Lagrangian
manifolds, the manifolds constructed of trajectories of 
material system elements (particles), which are obtained while describing
evolutionary processes in material media.

Below it will be shown that the properties of skew-symmetric forms 
depend on the properties of metric forms of the manifold on which these 
skew-symmetric forms are defined. The difference between exterior and 
evolutionary forms depends on the properties of the manifold's metric 
forms.

\subsection{Properties of the evolutionary differential
forms} 

As pointed out above, Cartan introduced the term ``exterior differential
forms" to denote skew-symmetric differential forms  (differential forms
with exterior multiplication). It was assumed that  these differential
forms are defined on manifolds which locally admit  a one-to-one mapping
into Euclidean subspaces [16] or into other  manifolds or submanifolds of
the same dimension. This means that these differential forms are defined
on manifolds with \emph{ closed metric forms}. In general, the theory of
exterior differential forms has been developed for differential forms
defined on such manifolds.

Skew-symmetric differential forms defined on manifolds with  \emph{
unclosed metric forms} have their own specific features that are  beyond
the scope of the modern theory of exterior differential forms. 
Therefore, it makes sense to introduce new terminology. In the  present
work we shall call them ``evolutionary differential forms".

The term ``evolutionary" is derived from the fact that differential 
forms defined on manifolds with unclosed metric forms possess 
the evolutionary properties. They are obtained from differential 
equations that describe any processes.
The coefficients of these forms depend on the evolutionary
variable, and the basis varies according to variation of  the form
coefficients.

Thus, exterior differential forms are defined on manifolds,  submanifolds
or on structures with closed metric forms. The evolutionary differential
forms are forms defined on manifolds with metric forms that are
unclosed and vary according to changing of the evolutionary variable.

[For the evolutionary differential forms we shall use a notation
with Greek indices, but for exterior differential forms we shall use
Latin indices as we have previously.]

Let us point cut some properties of evolutionary forms and show what 
their difference with  exterior differential forms consist in.

\subsection*{Specific features of the evolutionary forms 
differential }

An evolutionary form  differential of degree $p$ ($p$-form), 
as well as an exterior differential form, can be written as
$$
\omega^p=\sum_{\alpha_1\dots\alpha_p}a_{\alpha_1\dots\alpha_p}dx^{\alpha_1}\wedge
dx^{\alpha_2}\wedge\dots \wedge dx^{\alpha_p}\quad 0\leq p\leq n\eqno(2.1)
$$
where the local basis obeys the skew-symmetric condition 
$$
\begin{array}{l}
dx^{\alpha}\wedge dx^{\alpha}=0\\
dx^{\alpha}\wedge dx^{\beta}=-dx^{\beta}\wedge dx^{\alpha}\quad
\alpha\ne \beta
\end{array}
$$
(summation on repeated indices is implied).

But, the differential of the evolutionary form cannot be written in a 
manner similar to that described above for exterior differential forms 
(see formula (1.3)). In an evoltionary form differential there appears 
an additional term connected with the fact that the basis of the form
changes. An evoltionary differential form can be written as
$$
d\omega^p{=}\!\sum_{\alpha_1\dots\alpha_p}\!da_{\alpha_1\dots\alpha_p}dx^{\alpha_1}dx^{\alpha_2}\dots
dx^{\alpha_p}{+}\!\sum_{\alpha_1\dots\alpha_p}\!a_{\alpha_1\dots\alpha_p}d(dx^{\alpha_1}dx^{\alpha_2}\dots
dx^{\alpha_p})\eqno(2.2)
$$
where the second term is connected with the differential of the basis. 

(As is known, an exterior form differential does not include  the
differential of the basis. That is, in this case the differential of the
basis is equal to zero.)

[Hereinafter a symbol of summing $\sum$ and a symbol of exterior
multiplication $\wedge$ will be omitted. Summation over repeated indices
is implied.]

The properties of a skew-symmetric form differential depend on  the 
properties and specific features of the manifold on which skew-symmetric
differential form can be defined.  

The second term in the expression for the differential of skew-symmetric 
form connected with the differential of the basis is expressed in
terms of the metric form commutator. 
For differential forms defined on a manifold with unclosed metric
form, one has $d(dx^{\alpha_1}dx^{\alpha_2}\dots
dx^{\alpha_p})\neq 0$.
And for a manifold with a closed metric form, the following
$d(dx^{\alpha_1}dx^{\alpha_2}\dots dx^{\alpha_p}) = 0$ is valid.

That is, for differential forms defined on a manifold with unclosed 
metric form, the second term is nonzero, whereas for differential forms
defined on the manifold with closed metric form the second term vanishes.

For example, let us consider the first-degree form $\omega=a_\alpha
dx^\alpha$. The differential of this form can be written as
$$d\omega=K_{\alpha\beta}dx^\alpha dx^\beta$$
where $K_{\alpha\beta}=a_{\beta;\alpha}-a_{\alpha;\beta}$ are the
components of the commutator of the form $\omega$, and  
$a_{\beta;\alpha}$, $a_{\alpha;\beta}$ are
covariant derivatives. If we express the covariant derivatives in terms
of the connection coefficients (if it is possible), then they can be
written as $a_{\beta;\alpha}=\partial a_\beta/\partial
x^\alpha+\Gamma^\sigma_{\beta\alpha}a_\sigma$, where the first term
results from differentiating the form's coefficients, and the second term
results from differentiating the basis. (In the Euclidean space covariant
derivatives coincide with ordinary ones since, in this case, the
derivatives of the basis vanish). If we substitute the expressions for
covariant derivatives into the formula for commutator components, we
obtain the following expression for the commutator components of the form
$\omega$
$$
K_{\alpha\beta}=\left(\frac{\partial a_\beta}{\partial
x^\alpha}-\frac{\partial a_\alpha}{\partial
x^\beta}\right)+(\Gamma^\sigma_{\beta\alpha}-
\Gamma^\sigma_{\alpha\beta})a_\sigma\eqno(2.3)
$$
Here the expressions
$(\Gamma^\sigma_{\beta\alpha}-\Gamma^\sigma_{\alpha\beta})$ which entered
into the second term are just the components of the commutator of the
first-degree metric form. (For a nonintegrable manifold the connection
coefficients are nonsymmetric and the expressions
$(\Gamma^\sigma_{\beta\alpha}-\Gamma^\sigma_{\alpha\beta})$  are not equal
to zero).

That is, the corresponding metric form commutator will enter into the
differential form commutator.

If we substitute the expressions (2.3) for the skew-symmetric 
differential form commutator into formula for $d\omega$, we obtain the 
following expression for the differential of the first degree 
skew-symmetric form
$$
d\omega=\left(\frac{\partial a_\beta}{\partial
x^\alpha}-\frac{\partial a_\alpha}{\partial
x^\beta}\right)dx^\alpha dx^\beta+\left((\Gamma^\sigma_{\beta\alpha}-
\Gamma^\sigma_{\alpha\beta})a_\sigma\right)dx^\alpha dx^\beta
$$
The second term in the expression for the differential of a
skew-symmetric  form is connected with the differential of the manifold
metric form,  which is expressed in terms of the metric form commutator.

While deriving formula (2.3) for the differential form commutator
connection coefficients of a special type were used. However, a similar
result  can be obtained by applying a connection  of arbitrary type, or
by  using another means of finding the differential of the base
coordinates. For differential forms of any degree the metric form
commutator of corresponding degree will be included in the commutator of
the  skew-symmetric differential form.

As it is known [2,5], the differential of an exterior differential form
involves only a single term. There is no second term. This indicates that
the metric form commutator vanishes. In other words, the manifold, on
which the \emph{ exterior} differential form is defined, has a closed
metric form.

The differential of the evolutionary differential form, which is defined 
on a manifold with
unclosed metric forms, will contain two terms: the first term depends on
the differential form coefficients and the other depends on the
differential characteristics of the manifold.

Thus, the differentials and, correspondingly, the commutators of 
exterior and evolutionary forms are of different types.  
(As it will be shown below, this is precisely what determines 
the  characteristic properties and peculiarities of evolutionary forms.) 

What does this lead to? 

\subsection*{Non closure  of the evolutionary differential forms}

Commutators of \emph{ exterior} differential forms (differential forms
defined on manifolds with closed metric forms) contain only one term
obtained  from the derivatives of the differential form's coefficients.
Such a  commutator may be equal to zero. That is, the differential of
the  exterior form may vanish. This means that the exterior  differential
form may be closed.

In contrast to this case, evolutionary differential forms cannot be
closed.

Since  metric forms of the manifold are unclosed, the second term  of the
differential form commutator involves the metric form  commutator that is
not equal to zero.  Hence, the second term of the evolutionary form
commutator will be  nonzero. In addition, the terms of this commutator
have a different nature.  Such terms cannot compensate one another. For
this reason,  the differential form commutator proves to be nonzero. And
this means that the differential form defined on the manifold  with an
unclosed metric form cannot be closed.

Evolutionary differential forms are defined on manifolds with  unclosed
metric forms, and therefore, the evolutionary forms turn out to be 
unclosed.

\subsection{Specific features of the mathematical apparatus of evolutionary 
differential forms. Generate closed external forms}

Since evolutionary differential forms are unclosed, their mathematical 
apparatus  would not seem to possess capbilities connected with the
algebraic, group, invariant  and other such properties of closed exterior
differential forms. However, the mathematical apparatus  of evolutionary
forms proves to be significantly wider than expected because of the fact
that evolutionary differential forms can generate closed exterior
differential forms. 
Thus we obtain the answer to the question of how the closed exterior forms 
arise.

Such capabilities of evolutionary forms are due to the fact that
the mathematical apparatus of evolutionary differential forms includes 
unconventional elements, such one as nonidentical relations, degenerate 
transformations, transition from nonintegrable manifold to integrable one. 

\subsection*{Nonidentical relations of evolutionary differential forms}

In Section 1 it was shown that identical relations lie at the
basis of the mathematical apparatus of exterior differential forms. 
In contrast to this, nonidentical relations lie at the basis of the
mathematical apparatus of evolutionary differential forms. 

{\footnotesize [Identical relations establish exact correspondence between  
quantities (or objects) involved in a relation. It is possible
in  the case in which the quantities involved in the relation are invariants, 
in other words, measurable quantities. In a nonidentical  relation one of the
quantities is unmeasurable.]}  

A relation of evolutionary differential forms appears nonidentical  
as it involves an unclosed evolutionary form that is not invariant. 

Nonidentical relations appear in descriptions of any 
processes. They may be written as
$$
d\phi=\theta^p\eqno(2.4)
$$
Here $\theta^p$ is the $p$-degree evolutionary form that is an unclosed 
nonintegrable form, $\phi$ is some form of degree $(p-1)$, and the
differential $d\phi$ is a closed form of degree $p$.

Relation (2.4) is just of the same form as the identical relation (1.14)
from the mathematical apparatus of exterior differential forms presented 
in Section 1. However, on the right-hand side of the identical relation 
(1.14) stands a closed form, whereas the form on the right-hand side  of
nonidentical relation (2.4) is an unclosed one. (As will be shown below,
identical relations are satisfied only on pseudostructures, see (1.14')).

Relation (2.4) is an evolutionary relation as it involves an evolutionary
form. 

On the left-hand side of relation (2.4) stands a form differential, i.e.
a closed form that is an invariant object. On the right-hand side  stands
an unclosed form that is not an invariant  object. Such a relation cannot
be identical.

One can come to relation (2.4) by means of analyzing the integrability 
of the partial differential equation. An equation is integrable  if it
can be reduced to the form $d\psi=dU$. However it appears that, if the
equation is not subjected to an additional condition (the integrability
condition), it is reduced to the form (2.4), where $\omega$ is an
unclosed form and it cannot be expressed as a differential.

An example of a nonidentical relations among differential forms is a
functional relation (1.16) constructed by derivatives of the 
differential equation  (see Subsection 1.6).  In this
case  one has $\phi=u$ and $\omega=\theta$.

Here we present another derivation of a nonidentical relation (2.4) 
that  clarifies the physical meaning of this relation.

Let us consider this in terms of an example with a first degree
differential  form. A differential of a function of more than one
variables can be  an example of the first degree form. In this case the
function itself is the exterior form of zero degree.  The state  function
that specifies the state of a material system can serve as an example of
such a function. When the physical processes in a material  system are
being described, the state function may be unknown; but, its derivatives
may be known. The values of the function's derivatives may be equal to
some expressions that are obtained from the description of a real physical
process. And, the goal is to find the state function. 

Assume that $\psi$ is the desired state function that depends on the 
variables $x^\alpha$, and also assume that its derivatives in various
directions are known and equal to the quantities $a_\alpha$, namely:
$$
\frac{\partial\psi}{\partial x^\alpha}=a_\alpha\eqno(2.5)
$$
Let us set up the differential expression $(\partial\psi/\partial
x^\alpha)dx^\alpha$ (again, summation over repeated indices is implied).
This differential expression is equal to
$$
\frac{\partial\psi}{\partial x^\alpha}dx^\alpha=a_\alpha dx^\alpha\eqno(2.6)
$$
Here the left side of the expression is a differential of the
function $d\psi$, and the right side is the  differential form
$\omega=a_\alpha dx^\alpha$. Relation (2.6) can be written as 
$$
d\psi=\omega\eqno(2.7)
$$

It is evident that relation (2.7) is of the same type as (2.4) under the
condition that the differential form degrees  are equal to 1 (here the
right side is a first degree form, and the left side is a  differential
of the function, i.e. of the zero-degree form, which is the first degree
form as well). 

Relation (2.7) is a nonidentical because the differential form $\omega $ 
is an unclosed differential form. The commutator of this form is nonzero 
since the expressions $a_\alpha $ for the derivatives 
$(\partial\psi/\partial x^\alpha)$ are nonconjugated quantities.  They
are obtained from the description  of an actual physical process and are
unmeasurable quantities.

(Here it should be emphasized that the nonidentity of relation (2.7)  does
not mean that the mathematical description of the function's  variation
is not sufficiently accurate. The nonidentity of the  relation means that
the function's derivatives, whose values correspond to the real values in
physical processes, cannot be  consistent. While seeking a state function 
it is commonly assumed that  its derivatives are conjugated quantities, 
that is, their mixed  derivatives
commutative. But for physical processes the  expressions for these
derivatives are usually obtained independently of one another. And they
appear to be unmeasurable quantities, and hence they are not conjugated.  
Similar arguments may be also presented for the evolutionary relation
(2.4).) 

\bigskip

The nonidentical relation, which is evolutionary ones since it
includes the evolutionary form, plays a regulating role in the
process of generation of closed (inexact) exterior forms.

The process of obtaining closed exterior forms from evolutionary forms can
only proceed under the degenerate transformation. 

Since evolutionary forms are unclosed, the differential of
evolutionary forms is nonzero. And the differential of closed
exterior forms is equal to zero. From this it follows that the
transition from the evolutionary form to the closed exterior form
is only possible under the degenerate transformation, namely,
under the transformation which does not conserve the differential.

The degenerate transformation can only take place under additional conditions.
Such conditions can be realized in the evolutionary process described.

The nonidentical evolutionary relation, which appears to be a selfvarying 
relation, allows to find the conditions of such degenerate transformation.

\subsection*{Selfvariation of the evolutionary nonidentical relation}

An evolutionary nonidentical relation is selfvarying, because, firstly,
it is a nonidentical, namely, it contains two objects one of which appears
to be unmeasurable, and, secondly, it is an evolutionary relation,
namely, the variation of any object of the relation in some process leads
to a variation of another object; and, in turn, the variation of the
latter leads to variation of the former. Since one of the objects is an 
unmeasurable quantity, the other cannot be compared with the first  one,
and hence, the process of mutual variation cannot stop.  This process is
governed by the evolutionary form commutator. 

Varying an evolutionary form's coefficients leads to varying the first 
term of the commutator (see (2.3)). In accordance  with this variation,
it varies the second term, that is, also the metric form of the manifold
varies. Since the metric form's commutators specifies the manifold's
differential characteristics that are connected  with manifold
deformation (for example, the commutator of the  zero degree metric form
specifies the bend, that of second degree  specifies various
types of rotation, that of the third degree specifies  the curvature),
then it points to a manifold deformation. This means  that it varies
the evolutionary form basis. This, in turn, leads to  variation of the
evolutionary form, and the process of intervariation of the evolutionary
form and the basis is repeated. The processes of variation of the
evolutionary form and the basis are governed by the  evolutionary form's
commutator, and it is realized in accord with the evolutionary relation.

As a result the evolutionary relation appears to be a selfvarying
relation.  Selfvariation of the evolutionary relation proceeds by
exchange between the evolutionary form coefficients and manifold
characteristics. (This is an exchange between physical quantities and
space-time  characteristics, namely, between the algebraic and
geometrical  characteristics).  

The process of an evolutionary relation selfvariation cannot come to an 
end. This is indicated by the fact that both the evolutionary form 
commutator and the evolutionary relation involve unmeasurable
quantities. 

Just at such selfvariation  of an evolutionary relation it
can be realized the condition of degenerate transformation under
which a closed (inexact) exterior form can be obtained from the
evolutionary form.

\subsection*{Degenerate transformations. Realizations of 
pseudostructures and closed exterior differential forms} 
 
A distinction of the evolutionary form from the closed exterior form
consists in the fact that the evolutionary differential form is defined
on a manifold with unclosed
metric forms, and the closed exterior form can be defined only on a manifold
with closed metric forms.

Hence, it follows that a closed 
exterior form can be obtained from the evolutionary form only under 
\emph{ degenerate transformation}, when a transition from the   
manifold with unclosed metric forms (whose differential is nonzero) to 
integrable manifold with closed metric forms (for which the differential 
is zero) takes place. 

One can see that  vanishing the differential of the manifold metric form
is a condition of degenerate transformation.

For this reason a transition from an evolutionary 
form to a closed exterior form proceeds only when  
the differential or commutator of the metric form becomes equal to zero.

This can take place only discretely rather than identically. The 
coefficients of the commutator, or differential of the manifold metric 
form, as well as coefficients of the commutator or the evolutionary 
form differential, have different natures. Therefore, they cannot 
identically compensate one another (they cannot vanish identically). 
However, the commutators or differentials of a metric form can vanish 
given certain combinations of their coefficients. Such combinations 
(the conditions for a degenerate transformation) may be 
realized under selfvarying the evolutionary relation.

Since the closed metric form describes a pseudostructure,
vanishing of the metric form commutator (the realization of conditions 
for a degenerate transformation) and emergence of 
a closed metric form (the dual form) points to the realization 
of a pseudostructure and arises a closed inexact exterior form. 

It is clear that the transition from evolutionary form to closed one
is only possible if the closed metric form is realized. 

(The vanishing of one term of an evolutionary form commutator, namely, the
metric form commutator, leads to the fact that the second term of the
commutator also vanishes. This is due to the fact that the terms of an
evolutionary form commutator correlate with one another. An evolutionary
form commutator, and, correspondingly, the differential, vanish on a
pseudostructure, and this means that there arises a closed inexact
exterior form.)

Thus, if the conditions for a degenerate transformation are realized, 
from an unclosed evolutionary form, which the differential is nonzero
$d\theta^p\ne 0 $,
one can obtain a differential form closed on a pseudostructure. The
differential of this form equals zero. That is, it is realized the transition: 

$d\theta^p\ne 0 \to $ (a degenerate transformation) $\to d_\pi \theta^p=0$,
$d_\pi{}^*\theta^p=0$

\bigskip
{\bf Conditions of a degenerate transformation.} 
 
The Cauchy-Riemann conditions, the characteristic relations, the 
canonical relations, the Bianchi identities and others are examples of 
the conditions of degenerate transformations (vanishing the dual form
differential). 

Conditions of degenerate transformation, that is, additional conditions, can
be realized (under selfvariation of a nonidentical relation), for example, 
if there appear any symmetries  of the evolutionary or dual form coefficients  
or their derivatives. 

(At describing material systems this can be caused by an availability
of any degrees of freedom of material system.)

Corresponding to the conditions of degenerate transformation there is a
requirement, that certain functional expressions become equal to zero. Such
functional expressions (as it was pointed above in Subsection 1.6) are  
Jacobians, determinants, the Poisson brackets, residues, and others.

Mathematically, a degenerate transformation is realized as a transition from 
one frame of reference to another (nonequivalent) frame of reference. This is a 
transition from the frame of reference connected with the manifold whose metric 
forms are unclosed to the frame of reference being connected with a
pseudostructure. The first frame of reference cannot be inertial or a 
locally-inertial frame. The evolutionary form and nonidentical evolutionary 
relation are defined in the noninertial frame of reference. But the thereby 
obtained closed exterior form and the identical relation are obtained with 
respect to the locally-inertial frame of reference.

\subsection*{Obtaining an identical relation from a nonidentical} 

In Section 1 it was shown that identical relations of closed 
exterior form lie at the basis of many branches of mathematics.  

As one can see, the identical relations of closed exterior form follow
from nonidentical relations of unclosed evolutionary form.

On the pseudostructure $\pi$, from evolutionary relation (2.4) it follows
the relation
$$
d_\pi\phi=\theta_\pi^p\eqno(2.8)
$$
which turns out to be an identical relation. Indeed, since the form
$\theta_\pi^p$ is closed, on the pseudostructure this form turns
out to be a differential of some differential form. In other words,
this form can be written as $\theta_\pi^p=d_\pi\omega$. Relation (2.8)
is now written as $
d_\pi\phi=d_\pi\omega$. 
There are differentials on the left and right sides of this relation. This
means that the relation is an identical.

Thus, it is evident that under a degenerate transformation an identical 
relation on a pseudostructure can be obtained from an evolutionary nonidentical 
relation. This is due to the realization of a closed metric form, and 
correspondingly, to the realization of closed exterior form, and this points  
to the emergence of a pseudostructure and a conservaed quantity on 
the pseudostructure. The pseudostructure with conservated quantity constitutes 
a differential-geometrical  structure,  
which, as it has been shown in Subsection 1.2, is an invariant structure.
(The properties of such a structure will be defined below).

Under degenerate transformation an evolutionary form differential vanishes only
\emph{ on a pseudostructure}. It is an interior 
differential. The total differential of the evolutionary form is nonzero. The
evolutionary form remains unclosed, and for this reason, the original relation,
which contains the evolutionary form, remains a nonidentical.

Under realization of additional new conditions, a new identical relation can
be obtained. As a result, the nonidentical evolutionary relation can generate
identical  relations.

(It can be shown that all identical relations in the theory of exterior
differential forms are obtained from nonidentical relations by applying a
degenerate transformation.)

\bigskip
{\bf Connection  between nondegenerate transformations of exterior  forms and
degenerate transformations of evolutionary forms.} 

In the theory of closed
exterior forms only nondegenerate transformations, which conserve the
differential, are used. Degenerate transformations of evolutionary forms are
transformations that do not conserve the differential. Nevertheless, these
transformations are mutually connected. Degenerate transformations execute a
transition from an original deformed manifold to pseudostructures. And the
nondegenerate transformations execute a transition from one pseudostructure to
another. As the result of degenerate transformation, a closed inexact exterior
form arises from an unclosed evolutionary form, and this points 
to the generation of a differential-geometrical structure. 
But under nondegenerate transformation a transition from one closed
form to another takes place, and this points to the transition from one
differential-geometrical structure to another. 

\subsection*{Integration of a nonidentical evolutionary relation}

Since a closed exterior form is a differential (interior or total) of a  form 
of less by one degree, this allows us to integrate the closed form  and 
transition to a form of degree less by one. Such transitions are  possible 
in identical relations that connect a closed form with  a differential.

In can be shown that an integration and transitions with reduction by one
degree are also possible for a nonidentical relation (i.e., with an
evolutionary unclosed form). But this is possible only on condition of a
degenerate transformation are available.

Under degenerate transformation from the nonidentical evolutionary 
relation one obtains a relation being identical on pseudostructure. 
Since the right-hand side of such a relation can be expressed in terms 
of differential (as well as the left-hand side), one obtains a relation 
that can be integrated, and as a result he obtains a relation with the 
differential forms of less by one degree. 

The relation obtained after integration proves to be nonidentical 
as well. 

The resulting nonidentical relation of degree $(p-1)$ (relation that 
includes the forms of the degree $(p-1)$) can be integrated once again 
if the corresponding degenerate transformation has been realized and 
the identical relation has been formed. 

By sequential integrating the evolutionary relation of degree $p$ (in 
the case of realization of the corresponding degenerate transformations 
and forming the identical relation), one can get closed (on the 
pseudostructure) exterior forms of degree $k$, where $k$ ranges 
from $p$ to $0$. 

In this case one can see that under such integration the closed (on the 
pseudostructure) exterior forms, which depend on two parameters, are 
obtained. These parameters are the degree of evolutionary form $p$ 
(in the evolutionary relation) and the degree of created closed 
forms $k$. 

In addition to these parameters, an additiona parameter arisess, namely, the 
dimension of the space.

It is known that to a closed exterior differential forms of degree $k$
there correspond a skew-symmetric tensors of rank $k$ and to the corresponding
dual forms there is a pseudotensors of rank $(N-k)$, where $N$ is the 
dimension of the space. Pseudostructures correspond to such tensors,  
but only on the space formed.

\subsection{Functional possibilities of evolutionary forms}

The main peculiarity of evolutionary forms, which may be the deciding 
factor for mathematics, is the fact that 
the evolutionary forms generate closed inexact exterior forms, whose 
invariant properties lie at the basis of practically all invariant 
mathematical and physical formalisms.

The mathematical apparatus of exterior and evolutionary forms, which
basis involves nonidentical relations and degenerate
transformations, can describe transitions from nonconjugate
operators to conjugate ones and generation of various structures. 
(None of mathematical formalisms contains such possibilities.)

\subsection*{Mechanism of realization of conjugated objects and operators.}

To the closed exterior forms it can be assigned conjugated operators,  
whereas to the evolutionary forms there correspond nonconjugated 
operators. 
The transition from the evolutionary form to the closed exterior form
is that from nonconjugated operators to conjugated ones. This is expressed
as the  transition from the nonzero differential (unclosed evolutionary 
form) to the differential that equals zero  (closed exterior form).

Here it should be emphasized that the properties of deforming manifolds 
and skew-symmetric differential forms on such manifolds, namely, 
evolutionary forms, play a principal role in the process of conjugating.

The process of conjugating includes the following points:

1) selfvariation of nonidentical relation, namely, mutual variations 
of the evolutionary form coefficients (which have
the algebraic nature) and of manifold characteristics (which have the
geometric nature) described by nonidentical evolutionary relation, and

2) realization of the degenerate transformation.

Hence one can see that the process of
conjugating is a mutual exchange between the quantities
of different nature and the degenerate transformation under additional
conditions. Here it should be pointed that the condition of degenerate
transformation (vanishing some functional expressions like, for example,
Jacobians, determinants and so on) may be realized spontaneously
while selfvarying the nonidentical relation if any symmetries appear.
It is possible if the system (that is described by this relation) possesses
any degrees of freedom. 

One can see that the process of realization of conjugated operators or 
objects is described by the nontraditional mathematical apparatus, 
namely, by nonidentical relations and degenerate transformations.

As it has been shown above, closed exterior forms appear in many
mathematical formalisms. Practically all conjugated objects are
explicitly or implicitly connected with closed exterior forms. And
yet it can be shown that closed exterior forms are generated by
evolutionary differential forms, which are skew-symmetric
differential forms defined on the nonintegrable deforming varying 
manifolds.

{\footnotesize [It would be noted some specific features of mathematics. 
One branch of mathematics deals with conjugated operators or objects 
(algebra, geometry, the theory of groups, differential geometry, 
the theory of complex variables, exterior differential forms, and so on), 
whereas another deals with nonconjugated objects (differential and 
integral calculus, differential equations, topology and so on).

There are intersections between them. If the conjugacy conditions are known, 
one can obtain conjugated operators by imposing the conjugacy conditions 
on nonconjugated operators. However, in this case the questions of
how the conjugacy conditions are realized, what is the cause of
their origination and how the process of conjugating develops do
not solve in any  branch of mathematics.

The evolutionary differential forms answer these questions. They show 
how the transition from nonconjugated operators to conjugated 
ones proceeds.]}                                

\subsection*{Realization of differential-geometrical structures.} 

The process of generation of closed inexact
exterior forms describes thereby the process of
origination of the differential-geometrical structure  
which is an invariant structure (I-Structure).  

Obtaining differential-geometrical structures is a process of
conjugating the objects. Such process is, firstly,
a mutual exchange between the quantities of different nature
(for example, between the algebraic and geometric quantities or between
the physical and spatial quantities), and, secondly, the establishment
of exact correspondence (conjugacy) of these objects. This process 
is described by selfvariation of nonidentical relation and 
degenerate transformation.

\bigskip
{\bf Characteristics of the differential-geometrical structures
realized.} 

Since the closed exterior and dual differential forms, which
correspond to I-Structure arisen, were obtained from
the nonidentical relation that involves the evolutionary form, it is 
evident that the characteristics of such structure have to be connected: 

a) with those of the evolutionary form and of the deforming manifold
on which this form is defined, 

b) with the values of commutators of the
evolutionary form and the manifold metric form, and 

c) with the conditions of degenerate transformation as well.

The condition of degenerate transformation corresponds to a realization 
of the closed metric (dual) form and defines the pseudostructure.

Vanishing the interior commutator of the evolutionary form 
(on pseudostructure) corresponds to a realization of the closed (inexact) 
exterior form and points to emergence of conserved (invariant) quantity.

When I-Structure originates, the value of the total 
commutator of the evolutionary form containing two terms is nonzero. 
These terms define the following characteristics of I-Structures:

a) the first term of evolutionary form commutator (which is composed of 
the derivatives of the evolutionary form coefficients) defines the value of 
the discrete change of conserved quantity, that is, the quantum,
which the quantity conserved on the pseudostructure undergoes at the 
transition from one pseudostructure to another;

b) the second term (which is composed of the derivatives of
coefficients of the metric form connected with the manifold)
specifies the characteristics of I-Structures, which fixes the
character of the initial manifold deformation taking place before 
I-Structures  arose. 
(This characteristics fixes the deformation of original manifold that
proceeded in the process of originating the differential-geometrical 
structure and was described by selfvariation of nonidentical relation).

The discrete (quantum) change of conserved quantity proceeds in the 
direction that is normal to the pseudostructure. (Jumps of the 
derivatives normal to the potential surfaces are examples of such 
changes.)

(Above it has been noted that the evolutionary form and the nonidentical
relation are obtained while describing the physical processes that 
proceed in material systems. For this reason it is evident that
the characteristics of I-Structure must also be connected with
the characteristics of the material system being described.) 

\bigskip
{\bf Classification of differential-geometrical structures realized.}

The closed forms that correspond to I-Structures are generated by the
evolutionary relation which includes the evolutionary form of $p$ degree.
Therefore, the structures originated can be classified by the parameter 
$p$.

The other parameter is the degree of closed forms $k$ generated by
the nonidentical evolutionary relation.

Thus, one can see that I-structures, to which there are assigned
the closed (on the pseudostructure) exterior forms, can depend on two
parameters. These parameters are the degree of evolutionary form $p$
(in the evolutionary relation) and the degree of created closed
forms $k$.

In addition to these parameters, another parameter appears,
namely, the dimension of space. If the evolutionary relation
generates the closed forms of degrees $k$, to them there 
are assigned the pseudostructures of dimensions $(N-k)$, 
where $N$ is the space dimension.

\subsection*{Forming pseudometric and metric manifolds.}

At this point it should be noted that at every
stage of the evolutionary process it is realized only one element of
pseudostructure, namely, a certain mini-pseudostructure.

While varying the evolutionary variable the mini-pseudostructures form
the pseudostructure.
(The example of mini-pseudostructure is element of wave front. 
The element of wave front  made up the pseudostructure at its motion.)

Manifolds with closed metric forms are formed by pseudostructures. They
are obtained from the deforming manifolds with unclosed metric forms. 
In this case the initial deforming manifold (on which the evolutionary 
form is defined) and the manifold with closed metric forms originated 
(on which the closed exterior form is defined) are different  
objects.

It takes place the transition from the initial (deforming) manifold
with unclosed metric form to the pseudostructure, namely, to the
manifold with closed metric forms created. Mathematically this 
transition (the degenerate transformation) proceeds as
\emph{ a transition from one frame of reference to another, nonequivalent,
frame of reference.}

The pseudostructures, on which the closed \emph{ inexact} forms are
defined, form the pseudomanifolds.

To the transition from pseudomanifolds to metric space it is assigned
the transition from closed \emph{ inexact} differential forms 
to \emph{ exact} exterior differential forms.

It was shown above that the evolutionary relation of degree $p$ can
generate (with using the degenerate transformations) closed forms
of degrees $0,...,p$.  While generating closed forms of sequential
degrees  $k=p, k=p-1,..., k=0$ the pseudostructures of dimensions
$(n+1-k)$ are obtained. As a result of transition to the exact closed
form of zero degree the metric structure of the dimension $n+1$ is
obtained. 

Sections of the cotangent bundles (Yang-Mills fields), cohomologies 
by de Rham, singular cohomologies, pseudo-Riemannian and 
pseudo-Euclidean spaces, and others are examples of the pseudostructures 
and spaces that are formed by pseudostructures. Euclidean and Riemannian 
spaces are examples of metric
manifolds that are obtained when changing to the exact forms.
Here it should be noted that the examples of pseudometric spaces are
potential surfaces (surfaces of a simple layer, a double layer and so
on). In these cases the type of potential surfaces is connected with
the above listed parameters.

Conserved quantities (closed exterior inexact forms) defined on
pseudomanifolds (closed dual forms) constitute some fields. 
(The physical fields are the examples of such fields.) The fields of 
conserved quantities are formed from closed exterior forms at the 
same time when the manifolds are created from the pseudostructures.

Since the closed metric form is dual with respect to some closed exterior
differential form, the metric forms cannot become closed by themselves,
independently of the exterior differential form. This proves that
the manifolds with closed metric forms are connected with the closed
exterior differential forms.
This indicates that the fields of conserved quantities are
formed from closed exterior forms at the same time when the manifolds are 
created from the pseudostructures. The specific feature of manifolds with 
closed metric forms that have been formed is that they can carry some 
information.

\bigskip
One can see that the evolutionary forms possess the properties, which 
enable one to describe the evolutionary processes, namely, the processes 
of generating the differential-geometrical structures and forming 
manifolds. In other mathematical formalisms there are no such 
possibilities that the mathematical apparatus of evolutionary and 
exterior skew-symmetrical forms possesses.

\subsection*{Summary}

It is shown that the skew-symmetric differential forms play an unique 
role in mathematics.

The invariant properties of exterior skew-symmetric differential forms 
lie at the basis of practically all invariant mathematical.

The unique role of evolutionary skew-symmetric differential forms, which
were outlined in present work, relates to the fact that they generate
the closed exterior forms possessing invariant properties.

Due their properties and peculiarities the closed exterior forms and
evolutionary forms enable one to see the internal connection between
various branches of mathematics.

Many foundations of the mathematical apparatus of evolutionary
forms may occur to be of great importance for development of
mathematics. The nonidentical relations,
degenerate transformations, transitions from nonidentical
relations to identical ones, transitions from one frame of
reference to another (nonequivalent) frame, the generation of
closed inexact exterior forms and invariant structures, formatting
fields and manifolds, the transitions between closed inexact
exterior differential forms and exact forms and other phenomena
may find many applications in such branches of mathematics as the
qualitative theory of differential and integral equations,
differential geometry and topology, the theory of functions, the
theory of series, the theory of numbers, and others.

The evolutionary skew-symmetric differential forms may become a new
branch in mathematics.  They possess the possibilities that are 
contained in none of mathematical formalisms.

In the following work it will be shown an unique role of 
skew-symmetric differential forms in mathematical physics and field 
theory. Such role of skew-symmetric differential forms is due to the fact 
that they describe the properties of conservation laws.

1. Cartan E., Lecons sur les Invariants Integraux. -Paris, Hermann, 1922.  

2. Encyclopedia of Mathematics. -Moscow, Sov.~Encyc., 1979 (in Russian).

3. Novikov S.~P., Fomenko A.~P., Elements of the differential geometry and 
topology. -Moscow, Nauka, 1987 (in Russian). 

4. Bott R., Tu L.~W., Differential Forms in Algebraic Topology. 
Springer, NY, 1982. 

5. Schutz B.~F., Geometrical Methods of Mathematical Physics. Cambrige 
University Press, Cambrige, 1982.

6. Wheeler J.~A., Neutrino, Gravitation and Geometry. Bologna, 1960. 

7. Cartan E., Les Systemes Differentials Exterieus ef Leurs Application 
Geometriques. -Paris, Hermann, 1945. 

8. Konopleva N.~P. and Popov V.~N., The gauge fields. Moscow, Atomizdat, 1980 
(in Russian).

9. Smirnov V.~I., A course of higher mathematics. -Moscow, 
Tech.~Theor.~Lit. 1957, V.~4 (in Russian).

10. Haywood R.~W., Equilibrium Thermodynamics. Wiley Inc. 1980. 

11. Fock V.~A., Theory of space, time, and gravitation. -Moscow, 
Tech.~Theor.~Lit., 1955 (in Russian). 

12. Sternberg S., Lectures on Differential Geometry. -Englewood Cliffd , N.J.:
Prentiice-Hall,1964. 

13. Finikov S.~P.,  Method of the Exterior Differential Forms by Cartan in 
the Differential Geometry. Moscow-Leningrad, 1948 (in Russian).

14. Vladimirov V.~S., Equations of the mathematical physics. -Moscow, 
Nauka, 1988 (in Russian). 

15. Tonnelat M.-A., Les principles de la theorie electromagnetique 
et la relativite. Masson, Paris, 1959.

16. Efimov N.~V. Exterior Differntial Forms in the Euclidian Space. 
Moscow State Univ., 1971.

\end{document}